\documentclass[12pt]{article}
\usepackage{graphicx}
\usepackage{makecell}
\usepackage{algorithm}
\usepackage{algpseudocode}
\usepackage{float}
\usepackage{array}
\usepackage[utf8]{inputenc}
\usepackage{caption,graphicx,amsmath}
\graphicspath{ {images/} }
\usepackage{tocbibind}
\usepackage[colorlinks=true,
linkcolor=black,
filecolor=black,      
urlcolor=blue,
citecolor=blue]{hyperref}
\usepackage{epsfig,amsthm,latexsym,amssymb,amsgen}
\usepackage{amsthm}
\usepackage{cleveref}
\usepackage{bbm}
\usepackage{mathrsfs}
\usepackage{tabularx}
\usepackage[square, numbers]{natbib}
\usepackage{dcolumn}
\usepackage{bm}
\usepackage{amsfonts}
\usepackage{fontenc}
\usepackage{physics}
\usepackage{geometry}
\usepackage{subcaption}
\usepackage[toc,page]{appendix}
\usepackage{multirow}
\usepackage{enumitem}
\theoremstyle{plain}
\newtheorem{thm}{Theorem}[subsection]
\newtheorem{lemma}{Lemma}[subsection]
\newtheorem{prop}{Proposition}[subsection]
\newtheorem{cor}{Corollary}[subsection]
\newtheorem{rem}{Remark}[subsection]
\newtheorem{note}{Note}[subsection]
\theoremstyle{definition}
\newtheorem{definition}{Definition}[subsection]
\newtheorem*{example*}{Example}
\makeatletter
\newcommand{\thmnumbering}{%
	\ifbool{inSubsection}
	{\thesubsection.\arabic{thm}}
	{\thesection.\arabic{thm}}%
}

\newbool{inSubsection}
\pretocmd{\subsection}{\booltrue{inSubsection}}{}{}
\pretocmd{\section}{\boolfalse{inSubsection}}{}{}

\@addtoreset{thm}{section}
\@addtoreset{thm}{subsection}
\makeatother

\makeatletter
\newcommand{\lemmanumbering}{%
	\ifbool{inSubsection}
	{\thesubsection.\arabic{lemma}}
	{\thesection.\arabic{lemma}}%
}

\@addtoreset{lemma}{section}
\@addtoreset{lemma}{subsection}
\makeatother

\makeatletter
\newcommand{\cornumbering}{%
	\ifbool{inSubsection}
	{\thesubsection.\arabic{cor}}
	{\thesection.\arabic{cor}}%
}

\@addtoreset{cor}{section}
\@addtoreset{cor}{subsection}
\makeatother

\makeatletter
\newcommand{\propnumbering}{%
	\ifbool{inSubsection}
	{\thesubsection.\arabic{prop}}
	{\thesection.\arabic{prop}}%
}

\@addtoreset{prop}{section}
\@addtoreset{prop}{subsection}
\makeatother

\makeatletter
\newcommand{\notenumbering}{%
	\ifbool{inSubsection}
	{\thesubsection.\arabic{note}}
	{\thesection.\arabic{note}}%
}

\@addtoreset{note}{section}
\@addtoreset{note}{subsection}
\makeatother

\makeatletter
\newcommand{\remnumbering}{%
	\ifbool{inSubsection}
	{\thesubsection.\arabic{rem}}
	{\thesection.\arabic{rem}}%
}

\@addtoreset{rem}{section}
\@addtoreset{rem}{subsection}
\makeatother

\makeatletter
\newcommand{\definitionnumbering}{%
	\ifbool{inSubsection}
	{\thesubsection.\arabic{definition}}
	{\thesection.\arabic{definition}}%
}

\@addtoreset{definition}{section}
\@addtoreset{definition}{subsection}
\makeatother

\providecommand{\keywords}[1]
{
	\small	
	\textbf{\textit{Keywords---}} #1
}
\title{\textbf{On the non-zero divisor graph of the Hamilton quaternions over $\mathbb Z_{2^n}$}}
\author{\normalsize Gopika Govind \footnote{\href{mailto:gopikagovind03@gmail.com}{gopikagovind03@gmail.com}} , Chithra A.V \footnote{\href{mailto:chithra@nitc.ac.in}{chithra@nitc.ac.in}} \hspace{0.1cm} and \hspace{0.1cm} Manibharathi T.M.S\\
	\small	Department of Mathematics, National Institute of Technology Calicut,\\
	\small Kerala, India-673601}
\date{}
\begin{document}
	\maketitle	
	\begin{abstract}
		Let $R$ be a ring with unity. The non-zero divisor graph of $R$, $\Phi(R)$, is the graph with vertex set $R\backslash \{0,1,-1\}$, and two vertices $x$ and $y$ are adjacent if and only if either $xy$ or $yx$ is non-zero. In this article we associate $\Phi(R)$ to the ring of Hamilton quaternions over $\mathbb Z_{2^n}$, $\mathbb H(\mathbb Z_{2^n})$. The detailed structure of the elements in $\mathbb H(\mathbb Z_{2^n})$ is presented, based on which various structural properties of the graph $\Phi(\mathbb H(\mathbb Z_{2^n}))$, such as connectedness, adjacency of vertices, traversability, and planarity, are studied. Furthermore, we derive bounds for clique number and chromatic number, and study the concepts of matching and domination number.
	\end{abstract}
	\section{Introduction}
	\paragraph{}The concept of associating a graph with a commutative ring was introduced by Beck in 1988 \cite{beck1988coloring}. Later, this was modified by Anderson and Livingston \cite{andersonlivingston} in 1999, and is known as the classical zero-divisor graph of a commutative ring $R$, denoted by $\Gamma(R)$. The definition of zero-divisor graph was further extended to non-commutative rings by Redmond \cite{redmond2002zero} in 2002. Based on these ideas, many researchers have defined various graphs associated to different rings \cite{anderson2021graphs}.
	\paragraph{}One such graph that was recently defined, is the non-zero divisor graph of a ring $R$, denoted by $\Phi(R)$. It was introduced in 2020 by Kadem \textit{et al.} \cite{kadem2020non}. For any ring $R$, the simple graph $\Phi(R)$ has vertex set $V(R)=R \backslash \{0,1,-1\}$, and two distinct vertices $x$ and $y$ are adjacent if and only if either $xy\neq0$ or $yx\neq0$. \vspace{-0.5cm}
	\paragraph{} Our study focuses on the structural properties of the non-zero divisor graph of the Hamilton quaternions over $\mathbb Z_{2^n}$, where for a natural number $m$, $\mathbb Z_m$ is the set of integers modulo $m$. \vspace{-0.5cm}
	\paragraph{} Section 2 comprises the definitions and results that provide essential background information to the reader. In section 3 we derive the structure of elements in $\mathbb H(\mathbb Z_{2^n})$. Section 4 provides a detailed study of the graph $\Phi(\mathbb H(\mathbb Z_{2^n}))$. The section includes results regarding connectedness of the graph, adjacency of the vertices, properties like minimum degree, maximum degree, girth, and other structural properties as well.
	
	\section{Preliminaries}
	\vspace{0.1cm}
	\paragraph{} Let $G$ be a graph with vertex set $V(G)$ and edge set $E(G)$. $G$ is said to be \textit{simple}, if it has no loops and parallel edges. If $G$ contains an $x-y$ path, for all pair of vertices $x$, $y \in V(G)$, then $G$ is \textit{connected}. The \textit{minimum degree} of the graph $G$, $\delta(G)$, is the minimum among the degrees of all vertices in $G$. Similarly, the \textit{maximum degree} of $G$, $\Delta(G)$, is the maximum among the degrees of all vertices in $G$. The \textit{diameter} of $G$, denoted by $diam(G)$, is the greatest distance between any two vertices of $G$. A \textit{complete graph} is a graph in which every vertex is adjacent to each other. For a graph $G$ that is not complete, the \textit{connectivity} $\kappa(G)$ is the cardinality of a minimum vertex-cut. The \textit{edge-connectivity} $\lambda (G)$ of a nontrivial graph $G$ is the cardinality of a minimum edge-cut. A connected graph that contains a closed trail which involves all the edges in the graph (Eulerian circuit) is called an \textit{Eulerian graph}. If the graph $G$ contains a Hamiltonian cycle, that is, a cycle containing all the vertices in $G$, then $G$ is said to be a \textit{Hamiltonian graph}. The graph $G$ is called \textit{planar} if $G$ can be drawn in the plane so that no two of its edges cross each other. The \textit{clique number} of $G$, denoted by $\omega(G)$ is the number of vertices in the largest complete subgraph of $G$. The smallest number of colors in a proper coloring of the graph $G$ is known as the \textit{chromatic number} $\chi(G)$ of $G$.
	 A set of edges (vertices) in $G$ is said to be \textit{independent} if no two edges (vertices) in the set are adjacent. An independent set of edges in $G$ is known as a \textit{matching}. The \textit{edge independence number} $\alpha'(G)$ of the graph $G$ is the cardinality of a maximum matching. If a graph of order $2m$ has a matching of cardinality $m$, then it is called a \textit{perfect matching}. A vertex and an incident edge are said to \textit{cover} each other. A set of edges of $G$ that covers all the vertices of $G$ is known as an \textit{edge cover}. The minimum cardinality of an edge cover of $G$ is the \textit{edge covering number} $\beta'(G)$. The \textit{independence number} $\alpha(G)$ is the maximum cardinality of an independent set of vertices in $G$. Analogous to the concept of edge cover, we have a set of vertices that covers all edges in $G$, known as the \textit{vertex cover}. The \textit{vertex covering number} $\beta(G)$ is the minimum cardinality of a vertex cover in $G$. A  vertex in $G$ is said to \textit{dominate} itself and its neighbors. A set $S$ of vertices of $G$ is a \textit{dominating set} if every vertex of $G$ is dominated by some vertex in $S$. The cardinality of a minimum dominating set is the \textit{domination number} of $G$, denoted by $\gamma(G)$. A \textit{total dominating set} $S_t$ of $G$  is a set of vertices such that every vertex of $G$ is adjacent to at least one vertex in $S_t$. The \textit{total domination number} $\gamma_t(G)$ is the minimum cardinality of a total dominating set. All these notations are from \cite{chartrand2013first}.
	\begin{thm}\label{T} \cite{chartrand2013first}
		For every graph $G$, $\kappa(G)\leq \lambda(G) \leq \delta(G)$.
	\end{thm}
	\begin{thm}\label{t}\cite{chartrand2013first}
		For every graph $G$, $\chi(G)\geq \omega(G)$.
	\end{thm}
	\begin{thm}\label{Bt}\cite{chartrand2013first}
		For every connected graph $G$ that is not an odd cycle or a complete graph, $\chi(G)\leq \Delta(G)$.
	\end{thm}
	\begin{definition}\cite{kadem2020non} Let $R$ be a ring. The \textit{\textbf{non-zero divisor graph}}, denoted by $\Phi(R)$, is a simple graph, having a vertex set $V(R)=R\backslash \{0,1,-1\}$. Two distinct vertices $x$ and $y$ are adjacent if and only if either $xy \neq 0$ or $yx \neq 0$.
	\end{definition}
	\begin{example*} Let $R=M_2(\mathbb Z_2)$ be the ring of all $2 \times 2$ matrices with entries from $\mathbb Z_2$. The vertex set of $\Phi(M_2(\mathbb Z_2))$ is $M_2(\mathbb Z_2) \backslash \displaystyle \left\{ \begin{bmatrix}
			0&0\\0&0
		\end{bmatrix}, \begin{bmatrix}
			1&0\\0&1
		\end{bmatrix} \right\}$, and the graph is as shown in \hyperref[fig1]{Figure 1}.
		\begin{figure}[H]\label{fig1}
			\begin{center}
				\includegraphics[width=8.5cm,height=7cm]{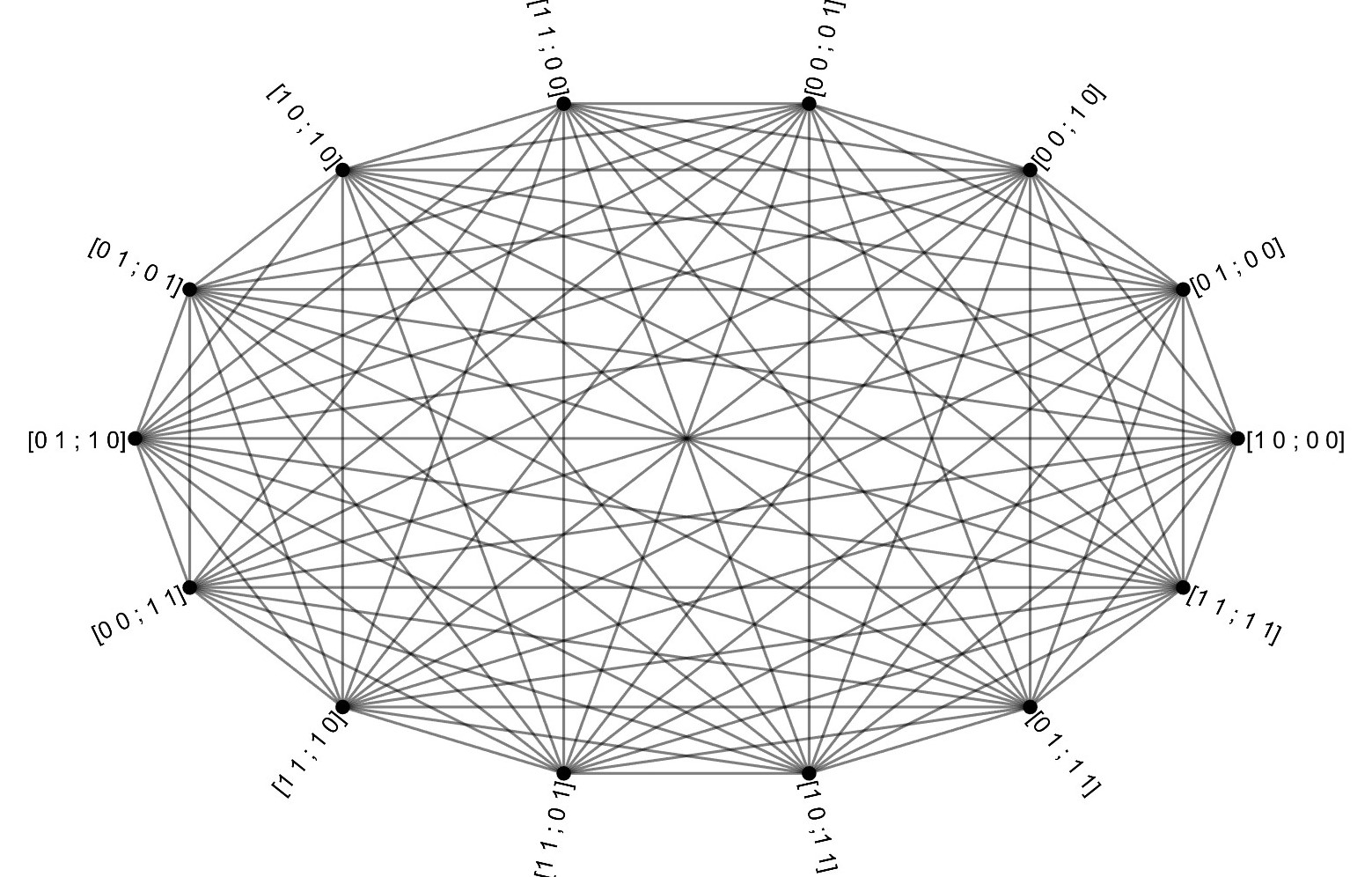}
			\end{center}
			\caption{\footnotesize Non-zero divisor graph of $M_2(\mathbb Z_2)$, $\Phi(M_2(\mathbb Z_2))$}
		\end{figure}
	\end{example*}
	\begin{definition}\cite{dummit2003third} \label{H(R)} Let $\mathbb H$ be the collection of the elements of the form $a+bi+cj+dk$, $a, b, c, d, a', b', c', d' \in \mathbb R$, where addition is defined component-wise by \vspace{-0.2cm}
		\begin{equation*}
			(a+bi+cj+dk)+(a'+b'i+c'j+d'k)=(a+a')+(b+b')i+(c+c')j+(d+d')k
			\vspace{-0.2cm}
		\end{equation*} 
		and multiplication is defined by expanding \vspace{-0.2cm}
		\begin{equation*}
			(a+bi+cj+dk)(a'+b'i+c'j+d'k) \vspace{-0.2cm}
		\end{equation*} 
		using the distributive law, and simplifying \vspace{-0.2cm}
		\begin{equation*}
			i^2=j^2=k^2=-1,\hspace{0.25cm} ij=-ji=k,\hspace{0.25cm} jk=-kj=i,\hspace{0.25cm} ki=-ik=j \vspace{-0.2cm}
		\end{equation*} 
		where the real number coefficients commute with $i$, $j$ and $k$.\\
		The collection $\mathbb H$ is known as the ring of \textit{\textbf{real Hamilton Quaternions}}.
		
	\end{definition}
	\begin{definition}\cite{dummit2003third} The \textbf{\textit{ring of \textbf{Hamilton Quaternions over integers modulo $m$}}, $\mathbb H(\mathbb Z_m)$}  can be defined by taking $a, b, c, d$ from $\mathbb Z_m$ in the definition of \hyperref[H(R)]{real Hamilton Quaternions}.\\
		The number of elements in $\mathbb H(\mathbb Z_m)$ is $m^4$.
	\end{definition}
	\begin{note}\label{n1}
		In $\mathbb H(\mathbb Z_{2^n})$, number of zero-divisors (including 0)= $2^{4n-1}$= number of units.
	\end{note}
	\begin{note}
		For the sake of brevity, throughout the article, we denote the element $a+bi+cj+dk \in \mathbb H(R)$ as $(a,b,c,d)$,  where $R$ is any ring.\\
		Consider the elements $(a_1,a_2,a_3,a_4)$, $(b_1,b_2,b_3,b_4) \in \mathbb H (R)$. Their product can be expressed by the following equations:
		\begin{equation}\label{eqns}
			\begin{split}
				&a_1b_1-a_2b_2-a_3b_3-a_4b_4=A\\
				&a_2b_1+a_1b_2-a_4b_3+a_3b_4=B\\
				&a_3b_1+a_4b_2+a_1b_3-a_2b_4=C\\
				&a_4b_1-a_3b_2+a_2b_3+a_1b_4=D
			\end{split}
		\end{equation}
		where each equation corresponds to the coefficient of $1,i,j,k$ in the product respectively.
	\end{note}
	\begin{rem}
		If $A,B,C$ and $D$ in {\normalfont(\ref{eqns})} are all 0, then 	$(a_1,a_2,a_3,a_4)$ and $(b_1,b_2,b_3,b_4)$ are zero-divisors in $\mathbb H (R)$.
	\end{rem}
	\begin{definition}
		The \textbf{\textit{non-zero divisor graph of Hamilton quaternions over $\mathbb Z_{2^n}$}}, denoted by $\Phi (\mathbb H (\mathbb Z_{2^n}))$, is the simple, undirected graph having $\mathbb H (\mathbb Z_{2^n}) \backslash \{(0,0,0,0),(1,0,0,0),(2^n-1,0,0,0)\}$ as the vertex set,  and two vertices $a=(a_1,a_2,a_3,a_4)$, $b=(b_1,b_2,b_3,b_4)$ are adjacent if $ab \neq 0$. \\
		Let $V(\Phi(\mathbb H(\mathbb Z_{2^n})))$ be the vertex set of $\Phi(\mathbb H(\mathbb Z_{2^n}))$. Then,\\
		\begin{equation*}
			|V(\Phi(\mathbb H(\mathbb Z_{2^n})))|=\begin{cases}
				2^{4n}-2, & n=1\\
				2^{4n}-3, & n>1 .
			\end{cases} 
		\end{equation*}
	\end{definition}
	\paragraph{}Throughout this article, the operations are taken to be modulo $2^n$, unless mentioned otherwise. 
	\begin{definition}\label{d1} \cite{NEWMAN1997367}
		For a commutative ring $R$ with unity, every matrix $A \in M_n(R)$ of rank $r$ is equivalent to a diagonal matrix $D$, given by
		\begin{equation*}
			D=diag(d_1,d_2,\cdots , d_r,0,\cdots,0),
		\end{equation*}
		where $d_i\neq 0$ for all $i \in \{1,2,\cdots,r\}$ and satisfy the divisibility sequence $d_1|d_2|\cdots |d_r$. The matrix $D$ is the \textbf{\textit{Smith normal form}} of $A$.\\
		It can be obtained by performing elementary row and column operations with respect to the operations in the ring $R$.
	\end{definition}
	
	\section{Structure of elements in $\mathbb H(\mathbb Z_{2^n})$}
	\paragraph{} In this section we obtain the structure of zero-divisors and units in $\mathbb H (\mathbb Z_{2^n})$. 
	\begin{thm}\label{t1}
		The elements in the ring $\mathbb H(\mathbb Z_{2^n})$ are as follows:
		\begin{itemize}
			\item[(i)] $(a_1,a_2,a_3,a_4)$, where all $a_i$'s are zero-divisors in $\mathbb Z_{2^n}$.
			\item[(ii)] $(b_1,b_2,b_3,b_4)$, where all $b_i$'s are units in $\mathbb Z_{2^n}$.
			\item[(iii)] All possible combinations of $2$ units and $2$ zero-divisors in $\mathbb Z_{2^n}$.
			\item[(iv)]  All possible combinations of $3$ units and $1$ zero-divisor in $\mathbb Z_{2^n}$. 
			\item[(v)] All possible combinations of $1$ unit and $3$ zero-divisors in $\mathbb Z_{2^n}$.
		\end{itemize}
	\end{thm}
	\begin{proof}
		\begin{itemize}
			\item[\textit{(i)}] Let $(a_1,a_2,a_3,a_4)\in \mathbb H(\mathbb Z_{2^n})$ such that $a_1$, $a_2$, $a_3$, $a_4$ are zero-divisors, thus 0 (not all simultaneously 0), or even, in $\mathbb Z_{2^n}$.\\
			Consider $(2^{n-1},2^{n-1},2^{n-1},2^{n-1}) \in \mathbb H(\mathbb Z_{2^n})$, we have,
			\begin{equation*}
				(2^{n-1},2^{n-1},2^{n-1},2^{n-1})(a_1,a_2,a_3,a_4)=(0,0,0,0). \vspace{0.1cm}
			\end{equation*} 
			Thus $(a_1,a_2,a_3,a_4)$ is a zero-divisor in $\mathbb H(\mathbb Z_{2^n})$.
			
			Similarly, we can see that the structures defined in \textit{(ii)} and \textit{(iii)} are also zero-divisors.
			\item[\textit{(iv)}] Consider an element $(a_1,a_2,a_3,a_4)$, where $a_1$, $a_2$, $a_3$ are units, and $a_4$ is a zero-divisor in $\mathbb Z_{2^n}$, and an arbitrary element $(b_1,b_2,b_3,b_4)$, in $\mathbb H(\mathbb Z_{2^n})$.\\
			Now, assume that $(a_1,a_2,a_3,a_4)(b_1,b_2,b_3,b_4)$ gives $A=B=C=D=0$ in (\ref{eqns})\\
			Thus (\ref{eqns}) becomes,
			\begin{equation}\label{2}
				a_1b_1-a_2b_2-a_3b_3-a_4b_4=0
			\end{equation}
			\begin{equation}\label{3}
				a_2b_1+a_1b_2-a_4b_3+a_3b_4=0
			\end{equation}
			\begin{equation}\label{4}
				a_3b_1+a_4b_2+a_1b_3-a_2b_4=0
			\end{equation}
			\begin{equation}\label{5}
				a_4b_1-a_3b_2+a_2b_3+a_1b_4=0 \vspace{0.2cm}
			\end{equation} 
			To prove that $(a_1,a_2,a_3,a_4)$ is a unit, it is enough to prove, that if \\ $(a_1,a_2,a_3,a_4)(b_1,b_2,b_3,b_4)=0$ then, $(b_1,b_2,b_3,b_4)=(0,0,0,0)$.\\
			We have,
			\begin{equation}\label{6}
				(a_1a_3+a_2a_4)b_1+(a_1a_4-a_2a_3)b_2-(a_3^2+a_4^2)b_3=0
			\end{equation}
			and
			\begin{equation}\label{7}
				(a_1a_3+a_2a_4)b_1+(a_1a_4-a_2a_3)b_2+(a_1^2+a_2^2)b_3=0
			\end{equation}
			respectively.\\
			From (\ref{6}) and (\ref{7}),
			\begin{equation}\label{8}
				(a_1^2+a_2^2+a_3^2+a_4^2)b_3=0. \vspace{0.1cm}
			\end{equation}
			Since $a_1^2,a_2^2$ and $a_3^2$ are odd, and $a_4^2$ is even (or 0), the sum $(a_1^2+a_2^2+a_3^2+a_4^2)$ in the LHS of (\ref{8}) is odd, and thus non-zero in $\mathbb Z_{2^n}$. Hence, for the RHS of (\ref{8}) to be 0, the only possibility is $b_3=0$.\\
			Using similar arguments, we can prove, that $b_1=b_2=b_4=0$.\\
			Proceeding in the same manner as above, we get, that the other possible structures defined in \textit{(iv)} and \textit{(v)} are also units.
		\end{itemize}
	\end{proof}
	\begin{rem}
		Since $\mathbb Z_{2^n}$ is a finite ring, correspondingly, $\mathbb H (\mathbb Z_{2^n})$ will also be finite. Thus all its elements are either a zero-divisor or a unit.
	\end{rem}
	\begin{example*}
		The number of elements in $\mathbb H(\mathbb Z_2)$ is $2^4=16$, out of which 8 are zero-divisors, and 8 are units. They are:\\
		Zero-divisors: $\{(0,0,0,0),(0,0,1,1),(0,1,0,1),(0,1,1,0),(1,0,0,1),(1,0,1,0),(1,1,0,0),(1,1,1,1)\}$.\\
		Units: $\{(1,0,0,0),(0,1,0,0),(0,0,1,0),(0,0,0,1),(0,1,1,1),(1,0,1,1),(1,1,0,1),(1,1,1,0)\}$.
	\end{example*}
	
	\section{The non-zero divisor graph of $\mathbb H(\mathbb Z_{2^n})$}
	\paragraph{} In this section, we study various structural properties of $\Phi(\mathbb H(\mathbb Z_{2^n}))$, like connectedness, diameter, degree of vertices, girth, planarity, traversibility, connectivity, chromatic number, matching, independence number, and domination number.\vspace{-0.5cm}
	\paragraph{}Throughout the section, we take two arbitrary vertices $a=(a_1,a_2,a_3,a_4)$ and $b=(b_1,b_2,b_3,b_4)$ of $\Phi(\mathbb H(\mathbb Z_{2^n}))$ as 
	$a_i=2^{k_i}u_i$, and $b_j=2^{l_j}v_j$, $0\leq k_i,l_j \leq n$ (not all simultaneously $n$), $i,j \in \{1,2,3,4\}$, where $u_i$'s and $v_j$'s are units in $\mathbb Z_{2^n}$, unless mentioned otherwise.
	\begin{thm}\cite{kadem2020non}\label{t2}
		Let $R$ be a ring. If $\Phi(R)$ is a connected non-zero divisor graph and $R\ncong \mathbb Z_2 \times \mathbb Z_4$, then $diam(\Phi(R)) \leq 2$.
	\end{thm}
	\begin{prop}\label{p1}
		Let \hspace{0.07cm} $\mathcal U(\mathbb H(\mathbb Z_{2^n}))=U(\mathbb H(\mathbb Z_{2^n})) \backslash \{(1,0,0,0),(2^n-1,0,0,0)\}$, where $U(\mathbb H(\mathbb Z_{2^n}))$ is the set of all units in $\mathbb H(\mathbb Z_{2^n})$. Then, all the vertices in $\mathcal U(\mathbb H(\mathbb Z_{2^n})) \subseteq V(\Phi(\mathbb H(\mathbb Z_{2^n})))$ is adjacent to every vertex in $\Phi(\mathbb H(\mathbb Z_{2^n}))$.
	\end{prop}
	\begin{rem}\label{r1}
		The vertices in \hspace{0.03cm} $\mathcal U(\mathbb H(\mathbb Z_{2^n}))$ induces a complete subgraph in $\Phi(\mathbb H(\mathbb Z_{2^n}))$, of order $2^{4n-1}-1$ when $n=1$, and $2^{4n-1}-2$ when $n>1$.
	\end{rem}
	\begin{prop}\label{p2}
		Vertices in $\Phi(\mathbb H(\mathbb Z_{2^n}))$, of the form $a=(a_1,a_2,a_3,a_4)$, $a_i=2^mu_i$, are not adjacent to vertices of the following forms:
		\begin{itemize}
			\item[(i)]  $b=(b_1,b_2,b_3,b_4)$, $b_j=2^{n-m}v_j$, $0 < m<n$, and
			\item[(ii)] $c=(c_1,c_2,c_3,c_4)$, $c_l
			=2^{n-(m+1)}w_l$, $0\leq m<n$,
		\end{itemize}
		where, $u_i$'s, $v_j$'s and $w_l$'s are units in $\mathbb Z_{2^n}$.
	\end{prop}
	\begin{prop}\label{p3}
		Vertices of the form $a=(a_1,a_2,a_3,a_4)$, where, $a_i=0$ or $2^{n-1}u_i$, $u_i$'s are units in $\mathbb Z_{2^n}$,  are not adjacent to any vertex of the same form. Furthermore, the vertex $(2^{n-1}u_1,2^{n-1}u_2,2^{n-1}u_3,2^{n-1}u_4)$ is not adjacent to any vertex that is a zero-divisor in $\Phi(\mathbb H(\mathbb Z_{2^n})))$.
	\end{prop}
	\begin{prop}\label{p4}
		For $n>2$, vertices of the form $a=(a_1,a_2,a_3,a_4)$ are not adjacent to any vertex of the form $b=(b_1,b_1,b_1,b_1)$, where, $a_i=2^{n-2}u_i$, $b_1=2^mv_1$, $0\leq m <n$, and $u_i$'s, $v_1$ are units in $\mathbb Z_{2^n}$.  
	\end{prop}
	\begin{cor}\label{c2}
		The graph $\Phi(\mathbb H(\mathbb Z_{2^n})))$ is not complete.
	\end{cor}
	\begin{thm}\label{t3}
		The graph $\Phi(\mathbb H(\mathbb Z_{2^n}))$ is connected with $diam(\Phi(\mathbb H(\mathbb Z_{2^n})))=2$, and $rad(\Phi(\mathbb H(\mathbb Z_{2^n})))=1$.
	\end{thm}
	\begin{proof}
		By \Cref{p1} and \Cref{t2}, we have, that $\Phi(\mathbb H(\mathbb Z_{2^n}))$ is connected with $diam(\Phi(\mathbb H(\mathbb Z_{2^n})))\leq 2$.\\
		Consider $c=(2^{n-1},2^{n-1},2^{n-1},2^{n-1})$ and another zero-divisor $a=(a_1,a_2,a_3,a_4)$ in $ V(\Phi(\mathbb H(\mathbb Z_{2^n})))$. By \Cref{p3}, we know that they are not adjacent.\\
		But, by \Cref{p1}, both these vertices are adjacent to all the units in $ V(\Phi(\mathbb H(\mathbb Z_{2^n})))$.\\
		Consider one such unit, say $b=(b_1,b_2,b_3,b_4)$. \\
		Then, we have the path $c-b-a$, of length 2, which is one of the shortest $c-a$ paths. \\
		Now, since all the units are adjacent to all vertices, $rad(\Phi(\mathbb H(\mathbb Z_{2^n})))=1$.
	\end{proof}
	\paragraph{}The following theorem gives the minimum degree $\delta$, maximum degree $\Delta$, and the girth $gr$ of $\Phi(\mathbb H(\mathbb Z_{2^n}))$.
	\begin{thm}\label{t4}
		For $\Phi(\mathbb H(\mathbb Z_{2^n}))$, the following properties hold:
		\begin{itemize}
			\item[(i)]  $\delta(\Phi(\mathbb H(\mathbb Z_{2^n})))=\begin{cases} 2^{4n-1}-1,& n=1\\2^{4n-1}-2,& n>1. \end{cases}$
			\item[(ii)]  $\Delta(\Phi(\mathbb H(\mathbb Z_{2^n})))=\begin{cases} 2^{4n}-3,& n=1\\2^{4n}-4,& n>1. \end{cases}$
			\item[(iii)] $gr(\Phi(\mathbb H(\mathbb Z_{2^n})))=3$.
		\end{itemize}
	\end{thm}
	\begin{proof}
		\begin{itemize}
			\item[\textit{(i)}] From \Cref{p3}, $deg((2^{n-1},2^{n-1},2^{n-1},2^{n-1}))=\begin{cases}
				2^{4n-1}-1, &n=1\\
				2^{4n-1}-2, &n>1,
			\end{cases}$
			that is, the number of units in $V(\Phi(\mathbb H(\mathbb Z_{2^n})))$.\\
			Since every other vertex is adjacent to all units, it follows that for an arbitrary vertex $a$,
			$deg(a) \geq \begin{cases}
				2^{4n-1}-1, &n=1\\
				2^{4n-1}-2, &n>1.
			\end{cases}$\\
			Then, $\delta(\Phi(\mathbb H(\mathbb Z_{2^n})))=\begin{cases} 2^{4n-1}-1,& n=1\\2^{4n-1}-2,& n>1. \end{cases}$
		\end{itemize}
		\vspace{0.2cm}
		\textit{(ii)} and \textit{(iii)} are direct consequences of \Cref{p1}.
	\end{proof}
	\begin{example*}
		The following figure shows the non-zero divisor graph of $\mathbb H(\mathbb Z_2)$.
		\begin{figure}[H] \label{fig2}
			\begin{center}
				\includegraphics[width=8.5cm,height=7cm]{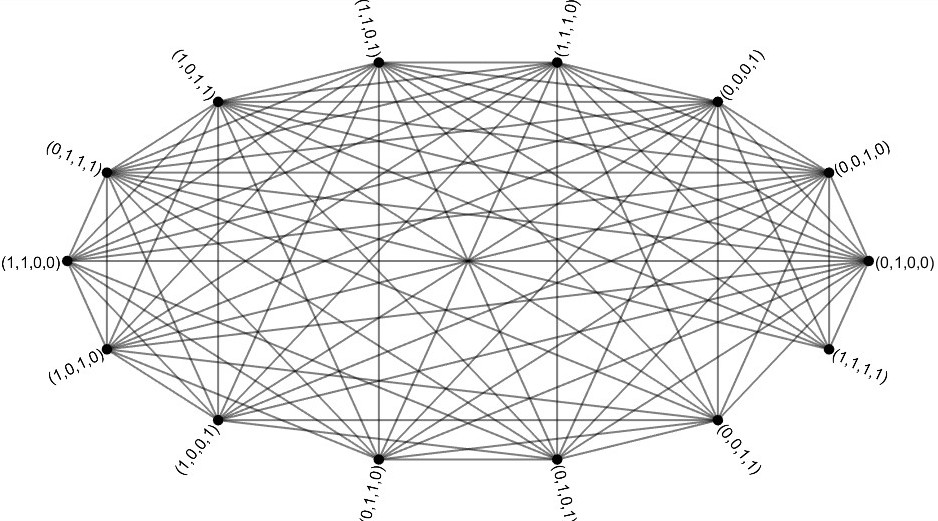}
			\end{center}
			\caption{\footnotesize{Non-zero divisor graph of $\mathbb H(\mathbb Z_2)$, $\Phi(\mathbb H(\mathbb Z_2))$}}
		\end{figure}
	\end{example*}
	\paragraph{} Let us now look into the structure-wise adjacency of certain vertices of $\Phi(\mathbb H(\mathbb Z_{2^n}))$.
	\begin{lemma}\label{l1}
		Suppose $a_i=2A_i$ denote zero-divisors, and $b_j=2B_j+1$ denote units in $Z_{2^n}$, respectively, where, $A_i, B_j \in Z_{2^n}$, with $i,j \in \{1,2,3,4\}$. For $n>1$, consider the vertices of the form:
		\begin{itemize}
			\item[(i)] $(a_1,a_2,b_1,b_2)$
			\item[(ii)] $(a_1,b_1,a_2,b_2)$
			\item[(iii)] $(a_1,b_1,b_2,a_2)$
			\item[(iv)] $(b_1,a_1,a_2,b_2)$
			\item[(v)] $(b_1,a_1,b_2,a_2)$
			\item[(vi)] $(b_1,b_2,a_1,a_2)$
			\item[(vii)] $(b_1,b_2,b_3,b_4)$
			\item[(viii)] $(2^lu,0,0,0)$, $(0,2^lu,0,0)$, $(0,0,2^lu,0)$, $(0,0,0,2^lu)$, where, $n>1$, $1\leq l<n-1$, and $u$ is a unit in $\mathbb Z_{2^n}$.
		\end{itemize}
		All the above vertices together induce a complete subgraph, of order $2^{n+1}(2^{3n-2}-2^{3n-5}+1)-8$, when $n>1$, in $\Phi(\mathbb H(\mathbb Z_{2^n}))$.
	\end{lemma}
	\begin{proof}
		Consider two vertices, one of the form $(a_1,a_2,b_1,b_2)$, and the other of the form $(a_3,b_3,a_4,b_4)$.\\
		Substituting these in (\ref{eqns}), we get at least one among $A,B,C,D$ as odd, and thus a non-zero output.\\
		This implies, that vertices of the form \textit{(i)} and \textit{(ii)} are adjacent in $\Phi(\mathbb H(\mathbb Z_{2^n}))$.\\
		Now, consider two vertices of the same form, say $(a_1,a_2,b_1,b_2)$ and $(a_3,a_4,b_3,b_4)$.\\
		Assume, that $A=B=C=D=0$ in (\ref{eqns}).\\
		Substituting the values of $a_i's$ and $b_j's$ in (\ref{eqns}), computing the first two equations modulo 4 and simplifying, we have
		\begin{equation*}
			\begin{split}
				B_1+B_2+B_3+B_4 &\equiv 1 \text{ } (mod \text{ } 2) \implies B_1+B_2+B_3+B_4=2p+1, \text{ for some $p\in \mathbb Z$},	\\
				B_1-B_2-B_3+B_4 &\equiv 0 \text{ } (mod \text{ } 2) \implies B_1-B_2-B_3+B_4=2q, \text{ for some $q\in \mathbb Z$}.
			\end{split}
		\end{equation*}
		From this we have, $2(B_2+B_3)=2(p-q)+1$, which is a contradiction.\\
		Thus, at least one among $A,B,C,D$ in (\ref{eqns}) should be non-zero.\\
		This implies, that vertices of the form \textit{(i)} are adjacent to each other.\\
		For $n>1$, using similar arguments we can obtain the other adjacencies as well.\\
		In the case of adjacency between vertices of the form \textit{(viii)}, substituting these structures pairwise in (\ref{eqns}), we get, that $(2^{l_1}u_1,0,0,0)$ is adjacent to $(2^{l_2}u_2,0,0,0)$ \\
		if and only if $l_1+l_2<n$. Same holds for other structures also. When $n>1$, there are $2^{n+1}-8$ such distinct vertices in $\Phi(\mathbb H(\mathbb Z_{2^n}))$.\\
		There are $2^{4n-1}$ zero-divisors in $\mathbb H(\mathbb Z_{2^n})$, out of which $2^{4(n-1)}-(2^{n+1}-8)$ zero-divisors of the form $(a_1,a_2,a_3,a_4)$ are not considered here.\\
		So, the order of the complete graph induced by the vertices of the forms \textit{(i)} to \textit{(viii)} $2^{n+1}(2^{3n-2}-2^{3n-5}+1)-8$, when $n>1$.
	\end{proof}
	\subsection{General formula for the adjacency and degree of vertices in $\Phi(\mathbb H(\mathbb Z_{2^n}))$}
	\paragraph{}In order to simplify the computations, and to obtain accurate results, for $n>1$, we reformulate (\ref{eqns}) as follows:
	\begin{equation}\label{refeqns}
		\begin{split}
			&\alpha_1\beta_1-\alpha_2\beta_2-\alpha_3\beta_3-\alpha_4\beta_4=A'\\
			&\alpha_2\beta_1+\alpha_1\beta_2-\alpha_4\beta_3+\alpha_3\beta_4=B'\\
			&\alpha_3\beta_1+\alpha_4\beta_2+\alpha_1\beta_3-\alpha_2\beta_4=C'\\
			&\alpha_4\beta_1-\alpha_3\beta_2+\alpha_2\beta_3+\alpha_1\beta_4=D'
		\end{split}
	\end{equation}
	where, 
	\begin{itemize}
		\item for $i,j \in \{1,2,3,4\}$, $\alpha_i=\displaystyle \frac{a_i}{2^k}$, $\beta_j=\displaystyle \frac{b_i}{2^l}$, with $k=\displaystyle \min_i\{k_i\}$ and $l=\displaystyle \min_j\{l_j\}$. 
		\item $A=2^{k+l}A'$, $B=2^{k+l}B'$, $C=2^{k+l}C'$ and $D=2^{k+l}D'$
	\end{itemize}
	\paragraph{}Throughout this article, $k$ and $l$ denote the values mentioned above, unless stated otherwise. Also, we denote the coefficient matrix of the above system as $\mathcal A$. The $2-$adic valuation of an integer $s$, which is the exponent of the highest power of 2 that divides $s$, is represented as $\nu_2(s)$, and we take
	$\nu=\min\{\nu_2(A'), \nu_2(B'),\nu_2(C'), \nu_2(D')\}$.
	\begin{thm}\label{t5}
		Consider two vertices in $\Phi(\mathbb H(\mathbb Z_{2^n}))$, say $a=(a_1,a_2,a_3,a_4)$ and $b=(b_1,b_2,b_3,b_4)$. Then, $a$ is adjacent to $b$ if and only if $k+l+\nu <n$.
	\end{thm}
	\begin{proof}
		Assume that $a$ is adjacent to $b$.\\
		This implies, that at least one among the outputs $A,B,C,D$ in (\ref{eqns}) will be non-zero.\\
		Without loss of generality, let $A\neq 0$. \\
		Then, $2^{k+l}A' \neq 0 \implies k+l+\nu_2(A')<n$.\\	
		Therefore, $k+l+\nu <n$.\\
		Conversely, assume that $k+l+\nu <n$.\\
		Without loss of generality, let $\nu=\nu_2(A')$.\\
		Then, $2^{k+l+\nu_2(A')}<2^n \implies2^{k+l}A' \neq 0$,
		which implies $A \neq 0$.\\
		Thus, $a$ is adjacent to $b$.
	\end{proof}
	\begin{rem}
		For any $n>1$, when $k=n-1$ and $l=0$, we have $\nu\geq 2$. Thus, $k+l<n$, and $k+l+\nu\geq n$, which implies, that though $k+l<n$, the corresponding vertices need not be adjacent. This emphasizes the relevance of the term $\nu$ in \Cref{t5}.
	\end{rem}
	\begin{thm}\label{t6}
		Let $a=(a_1,a_2,a_3,a_4)$ be an arbitrary vertex of $\Phi(\mathbb H(\mathbb Z_{2^n}))$, and $B=\{b\in(b_1,b_2,b_3,b_4)\in \mathbb H(\mathbb Z_{2^n}): ab=0 \}$. Then,
		\begin{itemize}
			\item[(i)] For $n=1$,
			\begin{equation*}
				deg(a)=\begin{cases}
					2^{4n}-|B|-1, & \text{if } a^2=0\\
					2^{4n}-|B|-2, & \text{if } a^2 \neq 0
				\end{cases}
			\end{equation*}
			\item[(ii)] For $n>1$,
			\begin{equation*}
				deg(a)=\begin{cases}
					2^{4n}-|B|-2, & \text{if } a^2=0\\
					2^{4n}-|B|-3, & \text{if } a^2 \neq 0
				\end{cases}
			\end{equation*}
		\end{itemize}
		where, $|B|=\displaystyle \prod_{i=1}^{4} gcd(d_i,2^n)$, in which $d_i$ is the $(i,i)^{th}$ entry in the Smith normal form of $A$.
	\end{thm}
	\begin{proof}
		We want the number of distinct elements in $V(\Phi(\mathbb H(\mathbb Z_{2^n})))$ which when multiplied with $a=(a_1,a_2,a_3,a_4)$ gives a non-zero result.\\
		For this, it is enough to subtract the number of distinct elements that produce zero output from the total number of vertices.\\
		Thus, 
		\begin{itemize}
			\item[\textit{(i)}] For $n=1$,
			\begin{equation*}
				deg(a)=\begin{cases}
					2^{4n}-|B|-1, & \text{if } a^2=0\\
					2^{4n}-|B|-2, & \text{if } a^2 \neq 0
				\end{cases}
			\end{equation*}
			\item[\textit{(ii)}] For $n>1$,
			\begin{equation*}
				deg(a)=\begin{cases}
					2^{4n}-|B|-2, & \text{if } a^2=0\\
					2^{4n}-|B|-3, & \text{if } a^2 \neq 0
				\end{cases}
			\end{equation*}
		\end{itemize}
		Now, we have to find $|B|$.\\
		Consider the Smith normal form of $A$, $diag(d_1,d_2,d_3,d_4)$ such that $d_i|d_{i+1}$.\\
		The number of solutions to the system of linear congruences $AX \equiv 0 \text{ }(mod \text{ }2^n)$ is \\
		$\displaystyle \prod_{i=1}^{4} gcd(d_i,2^n)$ (see section 5 in \cite{butson1955systems}).
	\end{proof}
	\paragraph{} Using a MATLAB program \cite{ferreira2009matlab}, we develop an algorithm to determine the neighbors and degree of a specific vertex. For a vertex $a$, the algorithm identifies all non-zero elements in $\mathbb H(\mathbb Z_{2^n})$ that yield a non-zero product when multiplied with $a$, except $(1,0,0,0)$, $(2^n-1,0,0,0)$, and $a$ itself, if $a^2\neq 0$. For values of $n\geq 5$, the program can be computationally expensive, due to the exponential growth of the number of vertices. This can be solved, by vectorizing or using parallel computing (parfor) to speed up the nested loops. The program can be summarized as in \Cref{alg:nonzero_product}.\vspace{0.2cm}
	\hrule \vspace{0.1cm}
		\captionof{algorithm}{Enumeration of quaternion tuples with non-zero product modulo $2^n$}
		\label{alg:nonzero_product} \vspace{-0.3cm} 
		\hrule \vspace{0.1cm}
		\begin{algorithmic}[1]
			\Require Quaternion tuple $a=(a_1,a_2,a_3,a_4)\in \mathbb Z_{2^n}^4$, positive integer $n$
			\Ensure Set $S$ of all tuples $b\in \mathbb Z_{2^n}^4$ such that $a\times b\neq (0,0,0,0)\pmod{2^n}$
			\State $m \gets 2^n$
			\State $S \gets \emptyset$
			\State $S \gets 0$
			\For{$b_1=0$ to $m-1$}
			\For{$b_2=0$ to $m-1$}
			\For{$b_3=0$ to $m-1$}
			\For{$b_4=0$ to $m-1$}
			\State $b \gets (b_1,b_2,b_3,b_4)$
			\State $A \gets (a_1b_1-a_2b_2-a_3b_3-a_4b_4)\bmod m$
			\State $B \gets (a_2b_1+a_1b_2-a_4b_3+a_3b_4)\bmod m$
			\State $C \gets (a_3b_1+a_4b_2+a_1b_3-a_2b_4)\bmod m$
			\State $D \gets (a_4b_1-a_3b_2+a_2b_3+a_1b_4)\bmod m$
			\If{$(A,B,C,D)\neq (0,0,0,0)$ \textbf{and} $b\neq a$ \textbf{and} $b\neq (1,0,0,0)$ \textbf{and} $b\neq (m-1,0,0,0)$}
			\State $S \gets S\cup\{b\}$
			\State $|S| \gets |S|+1$
			\EndIf
			\EndFor
			\EndFor
			\EndFor
			\EndFor
			\State \Return $(S, |S|)$
		\end{algorithmic}

	\subsection{Other structural properties of $\Phi(\mathbb H(\mathbb Z_{2^n}))$}
	\paragraph{} In this section we study about various structural properties of $\Phi(\mathbb H(\mathbb Z_{2^n}))$ such as connectivity, planarity, traversibility, matching etc.
	\begin{prop}\label{p5}
		When $n_1<n_2$, the graph $\Phi(\mathbb H(\mathbb Z_{2^{n_1}}))$ is a subgraph of $\Phi(\mathbb H(\mathbb Z_{2^{n_2}}))$.	
	\end{prop}
	\begin{prop}\label{p6}
		The graph $\Phi(\mathbb H(\mathbb Z_{2^n}))$ is not Eulerian.
	\end{prop}
	\begin{proof}
		The Smith normal form of $A$ corresponding to the vertex $a=(1,1,1,1)$ is $diag(1,0,0,0)$ when $n=1$, and $diag(1,2,2,4)$ when $n>1$. By \Cref{t6},
		\begin{equation*}
			deg(a)=\begin{cases}
				2^{4n}-2^3-1,& \text{ when } n=1\\
				2^{4n}-2^4-3,& \text{ when } n>1.
			\end{cases}
		\end{equation*}
		This indicates the existence of a vertex of odd degree in $\Phi(\mathbb H(\mathbb Z_{2^n}))$ for all $n$.
	\end{proof}
	\begin{thm}\label{t7}
		The graph $\Phi(\mathbb H(\mathbb Z_{2^n}))$ is Hamiltonian. 
	\end{thm}
	\begin{proof}
		We partition the vertex set of $\Phi(\mathbb H(\mathbb Z_{2^n}))$ into three disjoint sets $V_1,V_2$ and $V_3$, such that, the set $V_1$ contains all the non-zero zero-divisors in $\mathbb H(\mathbb Z_{2^n})$  of the form $(a_1,a_2,a_3,a_4)$, where $a_i's$ are all even in $\mathbb Z_{2^n}$. Therefore, $|V_1|=2^{4(n-1)}-1$.\\ $V_2$ contains all the other zero-divisors in $\mathbb H(\mathbb Z_{2^n})$. Then, $|V_2|=2^{4n-1}-2^{4(n-1)}$.\\ Finally, the set $V_3$ contains all the units in $\mathbb H(\mathbb Z_{2^n})$, except $(1,0,0,0)$ and $(2^n-1,0,0,0)$.
		This implies, $|V_3|=\begin{cases}
			2^{4n-1}-1,& \text{ when } n=1\\
			2^{4n-1}-2,& \text{ when } n>1.
		\end{cases}$
		\begin{itemize}
			\item[Case 1:] $n=1$.\\
			Note that $|V_1|=0$.\\
			By \Cref{p3}, the vertex $(1,1,1,1)$ is adjacent only to the units.\\
			From \hyperref[fig2]{Figure 2}, we can see, that each zero-divisor $z_i=(a_1,a_2,a_3,a_4)$ is adjacent to all other zero-divisors except $z_i'=(a_1,a_2,a_3,a_4)+(1,1,1,1)$, $i=1,2,3$.\\
			Let $u$ and $v$ be two units from the vertex set.
			We begin at the vertex $(1,1,1,1)$, move to $u$, then traverse the path $z_1-z_2-z_3-z_1'-z_2'-z_3'$, and then return to $(1,1,1,1)$ through the vertex $v$.\\
			This forms a Hamiltonian cycle.
			\item[Case 2:] $n>1$.\\
			Suppose $V_1=\{w_1,w_2,\cdots w_{2^{4(n-1)}-1}\}$,
			$V_2=\{v_1,v_2,\cdots,v_{2^{4n-1}-2^{4(n-1)}}\}$, and
			$V_3=\{u_1,u_2,\cdots,u_{2^{4n-1}-2}\}$.\\
			As a consequence of \Cref{l1}, we have the following:\\
			There exist a Hamiltonian path $H_1$ in the subgraph induced by the vertices in $V_2$, beginning from $v_1$ and ending at $v_{2^{4n-1}-2^{4(n-1)}}$. Furthermore, there is an edge $e_1=v_{2^{4n-1}-2^{4(n-1)}}-u_1$. An alternating path $P_1$ from $u_1$ to $w_{2^{4(n-1)}-1}$ is given by the sequence $u_1-w_1-u_2-w_2-\cdots - u_{2^{4n-1}-2}-w_{2^{4(n-1)}-1}$. Additionally, there exist an edge $e_2=w_{2^{4(n-1)}-1}-u_{2^{4(n-1)}}$, as well as a Hamiltonian path $H_2$ from $u_{2^{4(n-1)}}$ to $u_{2^{4n-1}-2}$. Also, there is an edge $e_3=u_{2^{4n-1}-2}-v_1$.\\
			$H_1 \cup e_1 \cup P_1 \cup e_2 \cup H_2 \cup e_3$ is a Hamiltonian cycle in $\Phi(\mathbb H(\mathbb Z_{2^n}))$.
		\end{itemize}
	\end{proof}
	\paragraph{} We use the above mentioned partition of the vertex set of $\Phi(\mathbb H(\mathbb Z_{2^n}))$, whenever necessary.
	\begin{prop}\label{p7}
		The graph $\Phi(\mathbb H(\mathbb Z_{2^n}))$ is non-planar.
	\end{prop}
	\begin{proof}
		By \Cref{p1}, for all values of $n$, $K_5$ is a subgraph of $\Phi(\mathbb H(\mathbb Z_{2^n}))$. Thus the result holds.
	\end{proof}
	\paragraph{} The following proposition gives the vertex connectivity $\kappa$, and edge connectivity $\lambda$ of $\Phi(\mathbb H(\mathbb Z_{2^n}))$, and establishes their relation with the number of units in the vertex set. The result is a consequence of Propositions \ref{p1}, \ref{p2} and \Cref{T}.
	\begin{prop}\label{p8}
		$\kappa(\Phi(\mathbb H(\mathbb Z_{2^n})))=\lambda(\Phi(\mathbb H(\mathbb Z_{2^n})))=|\mathcal U (\mathbb H(\mathbb Z_{2^n}))|$.
	\end{prop}
As a result of \Cref{r1}, \Cref{p3} and \Cref{l1}, we obtain a bound for the clique number of $\Phi(\mathbb H(\mathbb Z_{2^n}))$.
\begin{prop}\label{p9}
	For the graph $\Phi(\mathbb H(\mathbb Z_{2^n}))$,
	\begin{equation*}
		\omega(\Phi(\mathbb H(\mathbb Z_{2^n})))=2^{4n-1}+1, \text{ when } n=1	
	\end{equation*}
	and 	\begin{equation*}
		\omega(\Phi(\mathbb H(\mathbb Z_{2^n}))) \geq 2^{4n}-2^{4(n-1)}+2^{n+1}-10, \text{ when } n>1.
	\end{equation*}
\end{prop}
From Theorems \ref{t}, \ref{Bt}, \ref{t4} and \Cref{p9}, we get a bound for the chromatic number of $\Phi(\mathbb H(\mathbb Z_{2^n}))$ .
\begin{cor}
	For $n=1$, the chromatic number $\chi(\Phi(\mathbb H(\mathbb Z_{2^n})))$ of $\Phi(\mathbb H(\mathbb Z_{2^n}))$, is $2^{4n}-6$, and for $n>1$,
	\begin{equation*}
		2^{4n}-2^{4(n-1)}+2^{n+1}-10\leq \text{ } \chi(\Phi(\mathbb H(\mathbb Z_{2^n})))\text{ }\leq 2^{4n}-4.
	\end{equation*}
\end{cor}
\paragraph{} The following results provide the matching number of $\Phi(\mathbb H(\mathbb Z_{2^n}))$, and other related and analogous parameters.
\begin{thm}
	For the graph $\Phi(\mathbb H(\mathbb Z_{2^n}))$, the edge independence number is given by
	\begin{equation*}
		 \alpha'(\Phi(\mathbb H(\mathbb Z_{2^n})))=|\mathcal U(\mathbb H(\mathbb Z_{2^n})|
	\end{equation*}
	for all $n$.
\end{thm}
\begin{proof}
	Since vertices in $\mathcal U(\mathbb H(\mathbb Z_{2^n})$ are adjacent to all the vertices in $\Phi(\mathbb H(\mathbb Z_{2^n}))$, this implies, $|\mathcal U(\mathbb H(\mathbb Z_{2^n})| \leq \alpha'(\Phi(\mathbb H(\mathbb Z_{2^n})))$.\\
	Consider a subset $M$ of the edge set $E(\Phi(\mathbb H(\mathbb Z_{2^n})))$ such that an edge $e_i \in E(\Phi(\mathbb H(\mathbb Z_{2^n})))$ belongs to $M$ if one end-point of $e_i$ is $u_i \in \mathcal U(\mathbb H(\mathbb Z_{2^n})$, and the other end-point is a zero-divisor $z_i \in V(\Phi(\mathbb H(\mathbb Z_{2^n})))$, $1 \leq i \leq 2^{4n-1}-1$, when $n=1$, and $1 \leq i \leq 2^{4n-1}-2$, when $n>1$. Thus, $|M|=|\mathcal U(\mathbb H(\mathbb Z_{2^n})|$.\\
	Evidently, this is a matching. which means, $|M| \leq \alpha'(\Phi(\mathbb H(\mathbb Z_{2^n})))$.\\
	Now, assume that for some $e \in E(\Phi(\mathbb H(\mathbb Z_{2^n})))$, $M_0=M\cup \{e\}$ is a maximum matching.
	\begin{itemize}
		\item[Case 1: ] $e$ is an edge between two units or an edge between a zero-divisor and a unit.\\
		This implies, that $M_0$ is not an independent set.
		\item[Case 2: ] $e$ is an edge between two zero-divisors.\\
	Since for $n>1$, $|V_1|+|V_2|=|V_3|+1$, it follows that at least one end-point of $e$ coincides with that of another edge in $M$.
	\end{itemize}
	Thus, $M_0$ is not a maximum matching.\\
	Consequently, $M$ is a maximum matching, and $ \alpha'(\Phi(\mathbb H(\mathbb Z_{2^n})))=|M|=|\mathcal U(\mathbb H(\mathbb Z_{2^n})|$.
\end{proof}
\begin{cor}
	The graph $\Phi(\mathbb H(\mathbb Z_2))$ has a perfect matching.
\end{cor}
\begin{prop}
	For all values of $n$, let $M\subseteq E(\Phi(\mathbb H(\mathbb Z_{2^n})))$ denote the set of edges such that an edge $e_i \in E(\Phi(\mathbb H(\mathbb Z_{2^n})))$ belongs to $M$ if one end-point of $e_i$ is $u_i \in \mathcal U(\mathbb H(\mathbb Z_{2^n})$, and the other end-point is a zero-divisor $z_i \in V(\Phi(\mathbb H(\mathbb Z_{2^n})))$. When $n=1$, the indices satisfy $1 \leq i \leq 2^{4n-1}-1$, and when $n>1$, $1 \leq i \leq 2^{4n-1}-2$. Then $M$ is a minimum edge cover, and the edge covering number is
	\begin{equation*}
		 \beta'(\Phi(\mathbb H(\mathbb Z_{2^n})))=2^{4n-1}-1.
	\end{equation*}
\end{prop}\vspace{0.2cm}
As a consequence of Propositions \ref{p2}, \ref{p3}, and Lemma \ref{l1}, we get the following two results.
\begin{prop}
	For the graph $\Phi(\mathbb H(\mathbb Z_{2^n}))$, independence number,
	\begin{equation*}
			\alpha(\Phi(\mathbb H(\mathbb Z_{2^n})))=3,\text{ for }n=2,3
	\end{equation*}
	and
	\begin{equation*}
		\alpha(\Phi(\mathbb H(\mathbb Z_{2^n}))) \leq |V_1 \cup V_2|=2^{4n-1}-1 \text{ for }n>2.
	\end{equation*}
	The set $M=\{(1,1,1,1),(a_1,a_2,a_3,a_4),(b_1,b_2,b_3,b_4)\}$, where $(a_1,a_2,a_3,a_4)$ is a zero-divisor in $V(\Phi(\mathbb H(\mathbb Z_2)))$ such that $(a_1,a_2,a_3,a_4)$is not adjacent to $(b_1,b_2,b_3,b_4)$, is a maximum independent set in $\Phi(\mathbb H(\mathbb Z_2))$. Similarly, a set containing the vertices $(2,2,2,2)$, $(1,1,1,1)$ (or any zero-divisor with all odd entries), and $(0,0,2,2)$ (or any zero-divisor which contains two zeros and two even entries), is a maximum independent set in $\Phi(\mathbb H(\mathbb Z_{2^2}))$.
\end{prop}
\begin{example*}
	The sets $M_1=\{(1,1,1,1),(1,0,1,0),(0,1,0,1)\}$ and $M_2=\{(1,3,1,1),(2,0,2,0),(2,2,2,2)\}$ are maximum independent sets in $\Phi(\mathbb H(\mathbb Z_2))$ and $\Phi(\mathbb H(\mathbb Z_{2^2}))$ respectively.
\end{example*}
\begin{prop}
 For the graph $\Phi(\mathbb H(\mathbb Z_{2^n}))$, we have the vertex covering number
	\begin{equation*}
		\beta(\Phi(\mathbb H(\mathbb Z_{2^n})))=\begin{cases}
			2^{4n}-7,& \text{ for } n=1\\
			2^{4n}-18,& \text{ for } n=2,
		\end{cases}
	\end{equation*}
	and
	\begin{equation*}
		\beta(\Phi(\mathbb H(\mathbb Z_{2^n})))\leq 2^{4n}-4, \text{ for } n>2.
	\end{equation*}
\end{prop}
\begin{proof}
	For $n=1,2$ the proof is trivial.\\
	For $n>2$, since $V_2\cup V_3$ induces a complete subgraph, to cover all the edges in the induced subgraph, a vertex cover must contain $|V_2|+|V_3|-1$ vertices. It must also contain vertices from $V_1$ to cover the edges adjacent to only the vertices in $V_1$.\\
	This implies, $\beta(\Phi(\mathbb H(\mathbb Z_{2^n})))\leq |V_1|+|V_2|+|V_3|-1 = 2^{4n}-4$. 
\end{proof}
\paragraph{} The following proposition regarding the concept of domination is a direct consequence of \Cref{p1}.
 \begin{prop}
 	For all $n$, a minimum dominating set of the graph $\Phi(\mathbb H(\mathbb Z_{2^n}))$ is a singleton set consisting of any one unit, and a minimum total dominating set is a set consisting of exactly two units. Consequently, $\gamma(\Phi(\mathbb H(\mathbb Z_{2^n})))=1$, and $\gamma_t(\Phi(\mathbb H(\mathbb Z_{2^n})))=2$.
 \end{prop}

\section{Conclusion}
\paragraph{} In this paper, the structural properties of the non-zero divisor graph of the Hamilton quaternions over $\mathbb Z_{2^n}$ are studied. For all $n$,  the diameter of $\Phi(\mathbb H(\mathbb Z_{2^n}))$ is exactly 2. A general formula for the adjacency and degree of vertices in the graph has been derived. Also , we obtained  that $\Phi(\mathbb H(\mathbb Z_{2^n}))$ is non-planar, not Eulerian, but Hamiltonian. Moreover, bounds on the clique and chromatic number of $\Phi(\mathbb H(\mathbb Z_{2^n}))$ are presented. The concepts of matching and domination number have also been studied. Furthermore, we have developed a MATLAB algorithm to determine the neighbors and degree of any specific vertex in $\Phi(\mathbb H(\mathbb Z_{2^n}))$.
\nocite{*}
\bibliographystyle{plain} 
\renewcommand{\bibname}{References}
\bibliography{ref_on_nonzero_divisor_of_HZ2n}
\end{document}